\documentclass[11pt]{article}

\usepackage{amssymb}
\usepackage{amsfonts}
\usepackage{eucal}
\usepackage{amsmath}
\newtheorem{theorem}{Theorem}[section]
\newtheorem{proposition}[theorem]{Proposition}
\newtheorem{definition}[theorem]{Definition}

\newtheorem{lemma}[theorem]{Lemma}
\newtheorem{remark}[theorem]{Remark}
\newtheorem{example}[theorem]{Example}

\newenvironment{proof}[1][Proof]{\textbf{#1.} }{\ \rule{0.5em}{0.5em}}
\input pictex.sty

 \usepackage{cancel}
 \usepackage{xcolor}
 
 \usepackage{amsmath,color,ulem}
\newcommand{\stkout}[1]{\ifmmode\text{\sout{\ensuremath{#1}}}\else\sout{#1}\fi}

\begin{document}

 \title{Discrete Gr\"{o}nwall inequalities for demimartingales}
 	\author{Milto Hadjikyriakou\footnote{School of Sciences, University of Central Lancashire, Cyprus campus, 12-14 University Avenue, Pyla, 7080 Larnaka, Cyprus (email:mhadjikyriakou@uclan.ac.uk).}~~ and B.L.S. Prakasa Rao \footnote{CR RAO Advanced Institute of Mathematics, Statistics and Computer Science, Hyderabad 500046, India (e-mail: blsprao@gmail.com).}}	
 \maketitle

\begin{abstract}
The aim of this work is to obtain  discrete versions of stochastic Gr\"{o}nwall inequalities involving demimartingale sequences. The results generalize the respective theorems for martingales provided by Kruse and Scheutzow (2018) and Hendy et al. (2022). Moreover, we present an application which provides an upper bound for the a priori estimate of the backward Euler-Maruyama numerical scheme.
\end{abstract}


\textbf{MSC 2010:}

\section{Introduction}

The concept of positive association was introduced by Esary et al. in \cite{esary1967association} and it has been studied extensively due to its applicability in many different fields such as in physics, reliability, insurance mathematics, finance, biology etc. The definition is given below.

\begin{definition}
	A finite collection of random variables $X_1, \ldots, X_n$ is said to be (positively) associated if
	$$
	\operatorname{Cov}\left(f\left(X_1, \ldots, X_n\right), g\left(X_1, \ldots, X_n\right)\right) \geq 0
	$$
	for any componentwise nondecreasing functions $f, g$ on $\mathbb{R}^n$ such that the covariance is defined. An infinite collection is associated if every finite subcollection is associated.
\end{definition}

The concept of demimartingales introduced by Newman and Wright in \cite{NW1982} aims, among other purposes, to study the relation between this new dependence structure with sequences of partial sums of positively associated random variables and martingales.  The motivation for the definition of demimartingales was based on the following proposition which refers to mean zero positively associated random variables.

\begin{proposition}
	\label{prop1}Suppose $\left\{X_n, n \in \mathbb{N}\right\}$ are $L^1$, mean zero, associated random variables and $S_n=\sum_{i=1}^n X_i$. Then
	$$
	\mathbb{E}\left[\left(S_{j+1}-S_j\right) f\left(S_1, \ldots, S_j\right)\right] \geq 0, \, j=1,2, \ldots
	$$
	for all coordinatewise nondecreasing functions $f$.
\end{proposition}

\noindent The definition of demimartingales presented next is motivated by the previous proposition.

\begin{definition}
	A sequence of $L^1$ random variables $\left\{S_n, n \in \mathbb{N}\right\}$ is called a demimartingale if for all $j=1,2, \ldots$
	$$
	\mathbb{E}\left[\left(S_{j+1}-S_j\right) f\left(S_1, \ldots, S_j\right)\right] \geq 0
	$$
	for all componentwise nondecreasing functions $f$ whenever the expectation is defined. Moreover, if $f$ is assumed to be nonnegative, the sequence $\left\{S_n, n \in \mathbb{N}\right\}$ is called a demisubmartingale.
\end{definition}

It is clear by Proposition \ref{prop1} that the partial sum of mean zero associated random variables is a demimartingale. This conclusion generalizes, in some sense, the known result that the partial sums of mean zero independent random variables form a martingale sequence, since a martingale with the natural choice of $\sigma$-algebras is a demimartingale. This can easily be proven by using standard arguments, i.e.
$$
\begin{aligned}
\mathbb{E}\left[\left(S_{n+1}-S_n\right) f\left(S_1, \ldots, S_n\right)\right] & =\mathbb{E}\left\{E\left[\left(S_{n+1}-S_n\right) f\left(S_1, \ldots, S_n\right) \mid \mathcal{F}_n\right]\right\} \\
& =\mathbb{E}\left\{f\left(S_1, \ldots, S_n\right) \mathbb{E}\left[\left(S_{n+1}-S_n\right) \mid \mathcal{F}_n\right]\right\} \\
& =0
\end{aligned}
$$
where $\mathcal{F}_n=\sigma\left(X_1, \ldots, X_n\right)$. Similarly, it can be verified that a submartingale, with the natural choice of $\sigma$-algebras is also a demisubmartingale. The counterexample provided below, originally presented in \cite{MH2010}, proves that the converse statement is not always true. 

\begin{example}
	We define the random variables $\left\{X_1, X_2\right\}$ such that
	$$
	P\left(X_1=-1, X_2=-2\right)=p, P\left(X_1=1, X_2=2\right)=1-p
	$$
	where $0 \leq p \leq \frac{1}{2}$. Then $\left\{X_1, X_2\right\}$ is a demisubmartingale since for every nonnegative nondecreasing function $f$
	$$
	\mathbb{E}\left[\left(X_2-X_1\right) f\left(X_1\right)\right]=-p f(-1)+(1-p) f(1) \geq p(f(1)-f(-1))\geq 0. 
	$$
	Observe that $\left\{X_1, X_2\right\}$ is not a submartingale since
	$$
	\mathbb{E}\left[X_2 \mid X_1=-1\right]=\sum_{x_2=-2,2} x_2 P\left(X_2=x_2 \mid X_1=-1\right)=-2<-1.
	$$
\end{example}

It is important to highlight that the counterexample implies that the class of demimartingales is a wider class of random variables compared to martingales. Further, counterexamples provided in \cite{MH2010} also prove that not all demimartingales have positively associated demimartingale differences, a fact that proves that demimartingales form a class of random variables wider than the class of partial sums of zero mean positively associated random variables. Moreover, it worth to be mentioned that results obtained for demimartingales often generalize or even improve results available in the literature for (sub)martingales and positively associated random variables. Therefore this new class of random objects worth to be studied independently and in depth. For more on demimartingales, we refer the interested reader to the monograph of Prakasa Rao \cite{PR2012}.

\medskip

Recently, Kruse and Scheutzow (2018) in \cite{KS2018} and Hendy et. al. (2022) in \cite{HZS2022}, provided discrete versions of the stochastic Gr\"{o}nwall lemma involving a martingale. Considering the relation between martingales and demimartingales discussed above, the natural question that arises is whether similar results can be obtained for the class of demimartingales. The answer is provided in Section 2 of this paper where discrete versions of the Gr\"{o}nwall inequalities for demimartingales are obtained. In Section 3 we discuss an application to stochastic differential equations similar to the one presented in \cite{KS2018}.

\medskip

In what follows, the notation $x \wedge y$ represents the minimum between the real numbers $x$ and $y$, while for $d$-dimensional vectors $\mathbf{x}=(x_1, \ldots, x_d)$ and $\mathbf{y}=(y_1,\ldots, y_d)$ the notation $\langle \mathbf{x},\mathbf{y}\rangle$ stands for the inner product between vectors $\mathbf{x}$ and $\mathbf{y}$. As usual, for any real number $p$, the $p$-norm of a random variable $X$ is expressed as $\left\|X\right\|_p = \left(\mathbb{E}X^p \right)^{1/p}$. Throughout the paper we use the convention that $\sum_{j\in \emptyset} = 0$.

\medskip

The next result, provided here without a proof, is crucial for obtaining the main result of this work and can be found in \cite{C2003} and \cite{H2010} (check also Theorem 2.1.3 in \cite{PR2012}).

\begin{lemma}
	\label{stoptime}Let the sequence $\left\{S_n, n \geq 1\right\}$ be a demi(sub)martingale, $S_0=$ 0, and $\tau$ be a positive integer-valued random variable. Furthermore, suppose that the indicator function $I_{[\tau \leq j]}=h_j\left(S_1, \ldots, S_j\right)$ is a componentwise nonincreasing function of $S_1, \ldots, S_j$ for $j \geq 1$. Then the random sequence $\left\{S_j^*=S_{ \tau \wedge j}, j \geq 1\right\}$ is a demisubmartingale.
\end{lemma}

\section{Gr\"{o}nwall-type inequalities for demimartingales}
First, we provide a demimartingale inequality which provides an upper bound for the $p$-th moment of the supremum of a demimartingale in terms of the $p$-th moment of its negative infimum for $0<p<1$. The result generalizes Lemma 3 in \cite{KS2018} for the case of demimartingales and will be used for the proof of Gr\"{o}nwall inequalities.

\begin{lemma}
	\label{lemma1}Let $S_0, S_1, S_2, \ldots$ be a demimartingale sequence such that $S_0 \equiv 0$. Then for every $p\in (0,1)$ and every $n\in\mathbb{N}_0$ we have $$
	\mathbb{E}\left[\left(\sup _{0 \leq k \leq n} S_k\right)^p\right] \leq \frac{1}{1-p}\left(\mathbb{E}\left[-\inf _{0 \leq k \leq n} S_k\right]\right)^p.
	$$
\end{lemma}
\begin{proof}
	It can be easily proven that $\mathbb{E}S_j = \mathbb{E}S_{j+1},\, \forall j$. Therefore, in the case where $S_0 \equiv 0$ we have that $\mathbb{E}S_n   = 0 $ for all $n\in \mathbb{N}_0$. Therefore,
	$$
	0=\mathbb{E} S_n=\mathbb{E}\left(S_n \vee 0\right)-\mathbb{E}\left(\left(-S_n\right) \vee 0\right), \, \, \mbox{for all }\, \, n\in\mathbb{N}_0 
	$$
	which leads to
	$$
	\begin{aligned}
	\mathbb{E}\left(S_n \vee 0\right) & =\mathbb{E}\left(-S_n \vee 0\right) \leqslant \mathbb{E}\left(\sup_{0 \leqslant k \leqslant n}\left(-S_k\right)\right)=\mathbb{E}\left(- \inf_{0 \leqslant k \leqslant n}S_k\right).
	\end{aligned}
	$$
	Fix $x>0$ and define
	\[
	\tau_n = \max\{0\leq k\leq n: S_k\geq x\}
	\]
	for arbitrary chosen $n \in \mathbb{N}_0$. Set 
	$$S_j^*=S_{j \wedge \tau_n} \quad \forall j \geqslant 0.$$ Note that $$I\{\tau_n \geq j+1\}$$ is a componentwise nondecreasing function of $S_{1}, \ldots,S_j$  and by employing Lemma \ref{stoptime}, $S_j^*$ forms a demisubmartingale sequence. Thus, for the $n$ chosen above,
	\[
	\left\{\sup_{0\leq k\leq n} S_k \geq x\right\} = \left	\{\sup_{0\leq k\leq n} S_k^* \geq x\right\}.
	\]
	By utilizing the Doob's type inequality for a demisubmartingale sequence we have that for any $x>0$, 
	\begin{align*}
	&x P\left(\sup_{0\leq j\leq n}S_j\geq x\right)= x P\left(\sup_{0\leq j\leq n}S_j^*\geq x\right)\leq \mathbb{E}\left[S_n^*I\left\{\sup_{0\leq j\leq n}S_j^*\geq x\right\}\right] \leq \mathbb{E}(S_n^*\vee 0).
	\end{align*}
	As pointed out earlier, $\mathbb{E}S_j = 0 \, \,\forall j$, and hence $\mathbb{E}S_j^* = \mathbb{E}(S_{j\wedge \tau_n}) = 0, \, \forall j$ which leads to the conclusion that
	\begin{align*}
	&P\left(\sup _{0 \leqslant k \leqslant n} S_k \geqslant x\right) \leqslant  \frac{1}{x}\mathbb{E}\left(S_n^* \vee 0\right)\leqslant \frac{1}{x} \mathbb{E}\left(-\inf _{0 \leq j \leq n} S_j\right).
	\end{align*}
	Then, by applying standard arguments, we have that
	\begin{align*}
	&\mathbb{E}\left(\sup _{0 \leq j \leq 2} S_j\right)^p=\int_0^{\infty} P\left(\sup _{0 \leq j \leq n} S_j\geqslant x^{1/p}\right){\rm d}x \leq \int_0^{\infty}\left\{\frac{1}{x^{1/p}} \mathbb{E}\left(-\inf_{0 \leq j \leq n} S_j\right)\right\} \wedge 1 {\rm d}x \\
	& = \int_{0}^{Q^p}(Q x^{-1/p})\wedge 1{\rm d}x+ \int_{Q^p}^{\infty}(Q x^{-1/p})\wedge 1{\rm d}x,\,\,\mbox{where}\, \, Q = \mathbb{E}\left(-\inf_{0 \leq j \leq n} S_j\right)\\
	& =\frac{1}{1-p} \mathbb{E}\left( -\inf_{0 \leq j \leq n} S_j\right)^p.
	\end{align*}
	~~~~~~~~~~~~~~~~~~~~~~~~~~~~~~~~~~~~~~~~~~~~~~~~~~~~~~~~~~~~~~~~~~~~~~~~~~~~~~~~~~~~~~~~~~~~~~~~~~~~~~~~~~~~~~~~~~~
\end{proof}

\medskip 

\subsection{A discrete stochastic Gr\"{o}nwall inequality}
\noindent Next, we apply the demimartingale inequality provided in Lemma \ref{lemma1} to establish a first discrete Gr\"{o}nwall inequality.

\begin{theorem}
	\label{Gronw}Let $\left(X_n\right)_{n \in \mathbb{N}_0},\left(F_n\right)_{n \in \mathbb{N}_0}$, and $\left(G_n\right)_{n \in \mathbb{N}_0}$ be sequences of nonnegative random variables with $\mathbb{E}\left[X_0\right]<\infty$ such that
	\begin{equation}
	\label{cond}X_n \leq F_n+S_n+\sum_{k=0}^{n-1} G_k X_k, \quad \text { for all } n \in \mathbb{N}_0 
	\end{equation}
	where $\left(S_n\right)_{n \in \mathbb{N}_0}$ is a demimartingale such that $S_0 \equiv 0$. Then, for any $p \in(0,1)$ and $\mu, v \in[1, \infty]$ with $\frac{1}{\mu}+\frac{1}{v}=1$ and $p v<1$, it holds true that
	\begin{equation}
	\label{res1}\mathbb{E}\left[\sup _{0 \leq k \leq n} X_k^p\right] \leq\left(1+\frac{1}{1-v p}\right)^{\frac{1}{v}}\left\|\prod_{k=0}^{n-1}\left(1+G_k\right)^p\right\|_{\mu}\left(\mathbb{E}\left[\sup _{0 \leq k \leq n} F_k\right]\right)^p
	\end{equation}
	for all $n \in \mathbb{N}_0$. If $\left(G_n\right)_{n \in \mathbb{N}_0}$ is assumed to be a deterministic sequence of nonnegative real numbers, then for any $p \in(0,1)$ it holds true that
	\begin{equation}
	\label{res2}\mathbb{E}\left[\sup _{0 \leq k \leq n} X_k^p\right] \leq\left(1+\frac{1}{1-p}\right)\left(\prod_{k=0}^{n-1}\left(1+G_k\right)^p\right)\left(\mathbb{E}\left[\sup _{0 \leq k \leq n} F_k\right]\right)^p
	\end{equation}
	for all $n \in \mathbb{N}_0$.
\end{theorem}
\begin{proof}
	Let $\tilde{F}_n = \sup_{0\leq k\leq n}F_k$. Following the proof of Theorem 1 in \cite{KS2018} we have that 
	\[
	X_n \leq (\tilde{F}_n + L_n) \prod_{i=0}^{n-1}(1+G_i), \, \forall n\in \mathbb{N}_0
	\]
	where 
	\[
	L_n: = \sum_{k=0}^{n-1}(S_{k+1} -S_k)\prod_{j=0}^{k}(1+G_i)^{-1}. 
	\]
	Set $C_k  = \prod_{j=0}^{k}(1+G_i)^{-1}$. Then, 
	\begin{equation}
	\label{eq1}L_n = \sum_{k=0}^{n-1}C_k(S_{k+1}-S_{k})
	=C_{n-1}S_n + \sum_{k=1}^{n-1}(C_{k-1} - C_k)S_k.
	\end{equation}
	Observe that for any componentwise nondecreasing function $f$ we have that 
	\begin{align*}
	& \mathbb{E}\left[\left(L_{n+1}-L_n\right) f\left(L_1, \ldots, L_n\right)\right] = \mathbb{E}\left[C_{n}\left(S_{n+1}-S_n\right) f\left(L_1, \ldots, L_n\right)\right] = \mathbb{E}\left[\left(S_{n+1}-S_n\right) f_1\left(S_1, \ldots, S_n\right)\right]
	\end{align*}
	where $f_1\left(S_1, \ldots, S_n\right)=C_{n}f\left(L_1, \ldots, L_n\right)$. Since $G_k\geq0$ for all integers $k$
	$$
	C_{k+1}=\prod_{j=0}^{k+1}\left(1+G_j\right)^{-1}=\frac{C_k}{1+G_{k+1}}<C_k,
	$$
	due to \eqref{eq1} we  conclude that $f_1\left(S_1, \ldots, S_n\right)$ is a componentwise nondecreasing function of $S_1,\ldots,S_n$. Hence,  $(L_n)_{n\in \mathbb{N}_0}$ is a demimartingale. Let $\tilde{L}_n = \sup_{0\leq k\leq n}L_n$. Then, 
	\[
	\sup_{0 \leqslant k \leqslant n} X_k^p \leqslant\left(\tilde{F}_n+\tilde{L}_n\right)^p \prod_{i=0}^{n-1}\left(1+G_i\right)^p.
	\]
	By employing Holder's inequality we have that
	\begin{align*}
	& \mathbb{E}(\sup_{0 \leqslant k \leqslant n} X_k^p) \leq E\left[ \left(\tilde{F}_n+\tilde{L}_n\right)^p \prod_{i=0}^{n-1}\left(1+G_i\right)^p\right] \leqslant\left[\mathbb{E}\left(\prod_{i=0}^{n-1}\left(1+G_i\right)^p\right)^\mu\right]^{1 / \mu}\left[\mathbb{E}\left(\tilde{F}_n+\tilde{L}_n\right)^{v p}\right]^{1 / v}\\
	&\leq\left\|\prod_{i=0}^{n-1}\left(1+G_i\right)^p \right\|_\mu\left(\mathbb{E}(\tilde{F}_n^{vp})+\mathbb{E}(\tilde{L}_n^{vp})\right)^{\frac{1}{v}}.
	\end{align*}
	Since $(L_n)_{n\in \mathbb{N}_0}$ is a demimartingale, we can employ the result of Lemma \ref{lemma1} and get
	\begin{align*}
	&	\mathbb{E}(\tilde{L}_n^{vp}) \leq \frac{1}{1-vp} \mathbb{E}\left(-\inf _{0 \leq j \leq n} L_j\right)^{vp}  \leq \frac{1}{1-vp} \mathbb{E}(\tilde{F}_n^{vp}).
	\end{align*}
	The last inequality follows since due to the nonnegativeness of $X_n$ we have that $-L_n\leq \tilde{F}_n, \, \forall n\in \mathbb{N}_0$ which leads to $-\inf _{0 \leq j \leq n} L_j \leq \tilde{F}_n$. Thus,
	\begin{eqnarray*}
		\mathbb{E}\left(\sup_{0\leq k\leq n}X_k\right)^p &\leq& \left\|\prod_{i=0}^{n-1}\left(1+G_i\right)^p \right\|_\mu\left(\mathbb{E}(\tilde{F}_n^{vp})+\mathbb{E}(\tilde{L}_n^{vp})\right)^{\frac{1}{v}}\\
		&\leq& \left( 1+ \frac{1}{1-vp}\right)^{1/v}\left\|\prod_{i=0}^{n-1}\left(1+G_i\right)^p \right\|_\mu \left[ \mathbb{E}\left(\sup_{0\leq j\leq n}F_j\right)^{vp}\right]^{1/v}\\
		&\leq& \left( 1+ \frac{1}{1-vp}\right)^{1/v}\left\|\prod_{i=0}^{n-1}\left(1+G_i\right)^p \right\|_\mu \left[ \mathbb{E}\left(\sup_{0\leq j\leq n}F_j\right)^{p}\right]
	\end{eqnarray*}
	where the last inequality follows by applying Jensen's inequality. Inequality \eqref{res2} follows by applying the same steps for $\mu  = \infty$.  
\end{proof}

\begin{remark}
	It is important to highlight that, similar to the respective inequality for martingales, although the right hand side of \eqref{cond} depends on a demimartingale sequence, the upper bounds provided in \eqref{res1} and \eqref{res2} are uniform with respect to the demimartingale sequence.
\end{remark}

\subsection{A discrete fractional stochastic Gr\"{o}nwall inequality}
In this section we introduce a discrete fractional Gr\"{o}nwall inequality involving a sequence of mean zero associated random variables. This kind of inequalities are commonly used in the numerical analysis of multi-term stochastic time-fractional diffusion equations. The physical interpretation of the fractional derivative is that it represents a degree of memory in the diffusing material.

\medskip

Consider the following discretization for a time interval $[0,T]$: let $\tau$ be the temporal step-size, and let $N$ be a positive integer such that $\tau = T/N$. Define $t_n= n\tau$, for each $n = 0, 1, \ldots , N$.

\medskip

The uniform $L 1$-approximation for a multi-term Caputo temporal fractional derivative of orders $ 0<\beta_0 < \beta_1 < \cdots < \beta_m<1$ at the time $t_n$ is given by
$$
\begin{aligned}
\left.\sum_{r=0}^m q_r \frac{\partial^{\beta_r} f(t)}{\partial t^{\beta_r}}\right|_{t=t_n} =\sum_{r=0}^m q_r \frac{1}{\tau^{\beta_r} \Gamma\left(2-\beta_r\right)} \sum_{i=1}^n a_{n-i}^{\beta_r}\left(f\left(t_i\right)-f\left(t_{i-1}\right)\right)+r_\tau^n
\end{aligned}
$$
where  $a_j^{\beta_r}=(j+1)^{1-\beta_r}-j^{1-\beta_r}$ for each $j \geq 0$, $q_r$ are absolutely positive parameters  and $r_\tau^n$ is the truncation error (see for example \cite{HZS2022} and references therein). 

\begin{definition}
	Let $\left\{f_n\right\}_{n=0}^N$ be a sequence of real functions. The discrete time-fractional difference operator $D_\tau^{\beta_r}$ is given by
	$$
	D_\tau^{\beta_r} f_n=\frac{\tau^{-\beta_r}}{\Gamma\left(2-\beta_r\right)} \sum_{i=1}^n a_{n-i}^{\beta_r} \delta_t f_i=\frac{\tau^{-\beta_r}}{\Gamma\left(2-\beta_r\right)} \sum_{i=0}^n b_{n-i}^{\beta_r} f_i, \quad \forall n=1, \ldots, N
	$$
	where $\delta_t f_i=f_i-f_{i-1}$, and the constants are defined by $$b_0^{\beta_r}=a_0^{\beta_r}, b_n^{\beta_r}=-a_{n-1}^{\beta_r}, b_{n-i}^{\beta_r}=a_{n-i}^{\beta_r}-a_{n-i-1}^{\beta_r},$$ for each $i=1, \ldots, n-1$.
\end{definition}

\noindent In what follows, we denote
$$\lambda = \lambda_1 + \lambda_2/(2-2^{1-\beta_m})$$ 
where $\lambda_1$ and $\lambda_2$ are positive constants and let
$$W:=\sum_{r=0}^m q_r \frac{\tau^{1-\beta_r}}{\Gamma\left(2-\beta_r\right)} \sum_{j=1}^k a_{j-1}^{\beta_r}>0.$$
Finally, recall the Mittag-Leffler function which is of the form
\[
E_{\alpha}= \sum_{k=0}^{\infty}\dfrac{z^k}{\Gamma(1+k\alpha)}.
\]

\medskip

\noindent The result that follows provides a discrete version of a fractional Gr\"{o}nwall inequality involving a sequence of mean zero associated random variables. The particular result is motivated by the work presented in \cite{HZS2022} where a similar result involving a martingale sequence is presented (see Theorem 1).

\begin{theorem}
	Let $\left(Y_n\right)_{n \in \mathbb{N}}$ be a sequence of mean zero associated random variables and let  $\left(X_n\right)_{n \in \mathbb{N}},\left(F_n\right)_{n \in \mathbb{N}}$ be sequences of nonnegative random variables with $\mathbb{E}\left[X_0\right]<\infty$ such that
	\begin{equation}
	\label{cond_frac}\sum_{r=0}^m q_r D_\tau^{\beta_r} X_n \leq F_n+Y_n+\lambda_1 X_n+\lambda_2 X_{n-1}, \quad \forall n \geq 1
	\end{equation}
	where $q_r$ are positive integers for $r=0,1,\ldots,m$. Moreover, assume that $p \in(0,1)$ and $\mu, v \in[1, \infty]$ such that $\frac{1}{\mu}+\frac{1}{v}=1$ and $pv<1$. Then,
	$$
	\begin{aligned}
	\mathbb{E}\left[\sup _{1 \leq k \leq n} X_k^p\right] & \leq\left(1+\frac{1}{1-v p}\right)^{\frac{1}{v}}\left\|\left(2 E_{\beta_m}\left(2 \lambda t_n^{\beta_m} / q_m\right)\right)^p\right\|_{\mu} \\
	& \times\left(\mathbb{E}\left[\frac{\tau^{\beta_m}}{q_m \Gamma\left(1+\beta_m\right)} X_0 W\right]+\mathbb{E}\left[\frac{t_n^{\beta_m}}{q_m \Gamma\left(1+\beta_m\right)} \sup_{1\leq k\leq n}F_k\right]\right)^p.
	\end{aligned}
	$$
	If $\left(\lambda_q\right)_{q \in \mathbb{N}_0}$ is a deterministic sequence of nonnegative real numbers, then, for any $p \in(0,1)$,
	$$
	\begin{aligned}
	\mathbb{E}\left[\sup _{1 \leq k \leq n} X_k^p\right] & \leq\left(1+\frac{1}{1-p}\right)\left\|\left(2 E_{\beta_m}\left(2 \lambda t_n^{\beta_m} / q_m\right)\right)^p\right\|_{L^\mu(\Omega)} \\
	& \times\left(\mathbb{E}\left[\frac{\tau^{\beta_m}}{q_m \Gamma\left(1+\beta_m\right)} X_0 W\right]+\mathbb{E}\left[\frac{t_n^{\beta_m}}{q_m \Gamma\left(1+\beta_m\right)} \sup_{1\leq k\leq n}F_k\right]\right)^.
	\end{aligned}
	$$
\end{theorem}
\begin{proof}
	Let $\tilde{F}_n = \displaystyle\sup_{1 \leq k \leq n}F_j$. Note that since $(X_n)_{n\in\mathbb{N}}$ and $(F_n)_{n\in\mathbb{N}}$ are sequences of nonnegative random variables, we employ Lemma 3 in \cite{HZS2022} and for any integer $n$ we write $X_n$ as 
	\begin{align}
	\nonumber	X_n &\leq 2\left[\dfrac{\tau^{\beta_m}}{q_m\Gamma(1+\beta_m)}\left( \sum_{j=1}^{n}(F_j+X_j) + X_0W\right)  \right]E_{\beta_m}(2\lambda t_n^{\beta_m}/q_m)\\
	&\label{Xn_ineq}\leq 2\left[\dfrac{t_n^{\beta_m}}{q_m\Gamma(1+\beta_m)} \tilde{F}_n + S_n +\dfrac{\tau^{\beta_m}}{q_m\Gamma(1+\beta_m)}X_0W  \right]E_{\beta_m}(2\lambda t_n^{\beta_m}/q_m)
	\end{align}
	where $t_n = n\tau$ and $$S_n = \dfrac{\tau^{\beta_m}}{q_m\Gamma(1+\beta_m)}\sum_{j=1}^{n}X_j$$ forms a demimartingale sequence due to the assumption that $(X_n)_{n\in\mathbb{N}}$ are mean zero associated random variables. 
	Similar to the proof of Theorem \ref{Gronw}, we first apply Holder's inequality, i.e.
	\begin{align*}
	&	\mathbb{E}\left(\sup_{1\leq k\leq n}X_k^p\right)\leq \mathbb{E}\left[ \left(\dfrac{t_n^{\beta_m}}{q_m\Gamma(1+\beta_m)} \tilde{F}_n + S_n +\dfrac{\tau^{\beta_m}}{q_m\Gamma(1+\beta_m)}X_0W  \right)^p(2E_{\beta_m}(2\lambda t_n^{\beta_m}/q_m))^p\right]\\
	& \leq \left\| (2E_{\beta_m}(2\lambda t_n^{\beta_m}/q_m))^p  \right\|_\mu\left( \mathbb{E}\left( \dfrac{\tau^{\beta_m}}{q_m\Gamma(1+\beta_m)}X_0W\right)^{\nu p} +\mathbb{E}\left( \dfrac{t_n^{\beta_m}}{q_m\Gamma(1+\beta_m)} \tilde{F}_n\right)^{\nu p} + \mathbb{E}(\tilde{S}_n)^{\nu p} \right)^{\frac{1}{\nu}}
	\end{align*}
	where $\tilde{S}_n = \displaystyle\sup_{1 \leq k \leq n}S_j$. Observe that, since $X_n\geq 0$ for all $n$, then from \eqref{Xn_ineq} we have that 
	\[
	-S_n \leq \dfrac{t_n^{\beta_m}}{q_m\Gamma(1+\beta_m)} \tilde{F}_n  +\dfrac{\tau^{\beta_m}}{q_m\Gamma(1+\beta_m)}X_0W
	\]
	which leads to 
	\[
	-\inf_{1\leq k\leq n}S_k \leq \dfrac{t_n^{\beta_m}}{q_m\Gamma(1+\beta_m)} \tilde{F}_n  +\dfrac{\tau^{\beta_m}}{q_m\Gamma(1+\beta_m)}X_0W.
	\]
	Lemma \ref{lemma1} is applied to the term $\mathbb{E}(\tilde{S}_n)^{\nu p}$ and the desired result follows by applying similar steps as in the proof of Theorem 2.2.
\end{proof}

\begin{remark}
	The same result can be obtained if \eqref{cond_frac} is expressed as
	\[
	\sum_{r=0}^m q_r D_\tau^{\beta_r} X_n \leq F_n+T_n+\lambda_1 X_n+\lambda_2 X_{n-1}, \quad \forall n \geq 1
	\]
	where $T_n$ is a sequence of associated demimartingales. In this case, the expression $S_n$ is again a demimartingale sequence.
\end{remark}

\section{An application to numerical scheme for approximation of stochastic differential equations}
The Gr\"{o}nwall lemma is consider to be a result suitable for giving a priori bounds for the norm of
a solution of an inequality under the assumption of a linear growth estimate. As an application of the Gr\"{o}nwall inequality obtained in Section 2.1 we study the system of stochastic differential equations
$$\begin{aligned}
\mathrm{d} X(t) & =f(X(t)) \mathrm{d} t+g(X(t)) \mathrm{d} W(t), \quad t \in[0, T] \\
X(0) & =X_0,
\end{aligned}$$
where $f:\mathbb{R}^d\longrightarrow\mathbb{R}^d$ and $g:\mathbb{R}^{d}\longrightarrow\mathbb{R}^{d\times m}$ denote the drift and diffusion coefficient functions respectively.  Moreover, $W(t)$ denotes the column vector that consists of $m$ elements each one being a standard Wiener process defined from $[0,T]\times \Omega$ to $\mathbb{R}$. Following \cite{KS2018} we assume that the drift and the diffusion coefficient functions are related based on the coercivity condition described below:
\begin{equation}
\label{condf} \langle f(x),x\rangle+\dfrac{1}{2}|g(x)|^2\leq L(1+|x|^2),\, \forall x\in \mathbb{R}^{d}
\end{equation}
for $L\geq 0$, where $|\cdot|$ denotes the Euclidean norms on $\mathbb{R}^d$ and $\mathbb{R}^m$ as well as the Frobenius norm in case of matrices from $\mathbb{R}^{d\times m}$. Consider the numerical scheme given below known as the backward Euler–Maruyama method:
$$
\begin{aligned}
Y^{j+1} & =Y^j+h f\left(Y^{j+1}\right)+g(Y^j)\Delta_h W^{j+1}, \quad j=1, \ldots, N_h, \\
Y^0 & =X_0,
\end{aligned}
$$
where $h \in(0,1)$ denotes the equidistant step size and $N_h \in \mathbb{N}$ is determined by $N_h h \leq T<\left(N_h+1\right) h$. The stochastic increment is given by $\Delta_h W^{j+1}=W\left(t_{j+1}\right)-W\left(t_j\right)$, where $t_j=j h$. The Gr\"{o}nwall lemma will be used to obtain an upper bound for the $L^{2p}$-norm of the supremum of the sequence $(Y^j)_{0\leq j\leq N_h}$.

\begin{proposition}
	Let $h_0 \in\left(0,(2 L)^{-1}\right)$ denote an upper step size bound. For every $p \in(0,1)$ and for every $(\mathcal{F} _{nh})_{n\in \mathbb{N}_0}$ adapted process $\left(Y^n\right)_{n \in \mathbb{N}_0}$ satisfying the backward Euler–Maruyama scheme with $h \in\left(0, h_0\right)$ we have
	$$
	\begin{aligned}
	\left\|\sup _{0 \leq j \leq N_h}\left|Y^j\right|\right\|_{2 p} 
	\leq\left(\frac{2-p}{1-p}\right)^{\frac{1}{2 p}} \exp \left(\frac{L T}{1-2 h_0 L}\right)\left(\left|X_0\right|^2+\left(1-2 h_0 L\right)^{-1}\left(h_0\left|g\left(X_0\right)\right|^2+2 L T\right)\right)^{\frac{1}{2}}.
	\end{aligned}
	$$
\end{proposition}
\begin{proof}
	Following the same steps as in the proof of Proposition 5 in \cite{KS2018} we first obtain the following inequality for the sequence $\left(Y^n\right)_{n \in \mathbb{N}_0}$
	$$
	(1-2 h L)\left|Y^n\right|^2+h\left|g\left(Y^n\right)\right|^2 \leq(1-2 h L)\left|Y^0\right|^2+h\left|g\left(Y^0\right)\right|^2+2 L t_n+\sum_{j=0}^{n-1} Z^{j+1}+2 h L \sum_{j=0}^{n-1}\left|Y^j\right|^2
	$$
	where 
	$$Z^{j+1}:=\left|g\left(Y^j\right) \Delta_h W^{j+1}\right|^2-h\left|g\left(Y^j\right)\right|^2+2\left\langle g\left(Y^j\right) \Delta_h W^{j+1}, Y^j\right\rangle.$$
	The inequality provided above can be written in the form 
	$$
	X_n \leq F_n+S_n+\sum_{j=0}^{n-1} G_j X_j, \quad \text { for all } n \in \mathbb{N}_0
	$$
	where
	$$
	\begin{aligned}
	X_n & :=\left|Y^n\right|^2+h\left|g\left(Y^n\right)\right|^2 \\
	F_n & :=\left|Y^0\right|^2+\left(1-2 h_0 L\right)^{-1}\left(h_0\left|g\left(Y^0\right)\right|^2+2 L t_n\right) \\
	S_n & :=\left(1-2 h_0 L\right)^{-1} \sum_{j=0}^{n-1} Z^{j+1} \\
	G_n & :=\left(1-2 h_0 L\right)^{-1} 2 h L.
	\end{aligned}
	$$ 
	Since $0<(1-2h_0L)\leq(1-2hL)\leq 1$, sequences $(F_n)_{n\in\mathbb{N}_0}$ and $(G_n)_{n\in\mathbb{N}_0}$ satisfy the conditions of Theorem \ref{Gronw} and it remains to prove that $(S_n)_{n\in\mathbb{N}_0}$ forms a demimartingale. Note that the condition  $S_{0}\equiv 0$ is satisfied (recall the convention $\sum_{j\in \emptyset} = 0$). 
	
	\medskip
	
	\noindent It is given that $\left(Y^n\right)_{n \in N_0}$ is $\left(\mathcal{F}_{n h}\right)_{n \in N_0}$-adapted and recall that $t_n = nh$ where $h\in (0,1)$. Based on the discussion in the proof of Proposition 5 in \cite{KS2018}, $Z^{n}$ and, hence, $S_n$ are integrable random variables. Then, for any $f$ componentwise nondecreasing function, and by conditioning on $\mathcal{F}_{t_n}$, we have that
	$$
	\begin{aligned}
	& E\left[\left(S_{n+1}-S_n\right) f\left(S_1, \ldots, S_n\right)\right] 
	=\left(1-2 h_0 L\right)^{-1} E\left[Z^{n+1} f\left(S_1, S_2, \ldots, S_n\right)\right] \\
	& =\left(1-2 h_0 L\right)^{-1} E\left[Z^{n+1} f^*\left(Z^1, Z^2, \ldots, Z^n\right)\right] =\left(1-2 h_0 L\right)^{-1} E\left[E\left(Z^{n+1} f^*\left(Z^1, Z^2, \ldots, Z^n\right) \mid \mathcal{F}_{t_n}\right)\right] \\
	& =(1-2 h_0L)^{-1} E\left[f^*\left(Z, Z_1^2, \ldots, Z^n\right) E\left(Z^{n+1} \mid \mathcal{F}_{t_n}\right)\right]=0
	\end{aligned}
	$$
	since Kruse and Scheutzow in \cite{KS2018} proved that $E\left(Z^{n+1} \mid \mathcal{F}_{t_n}\right) = 0$ which concludes that the sequence $(S_n)_{n\in\mathbb{N}_0}$ satisfies the definition of a demimartingale sequence. 
	
	\medskip
	
	\noindent  Because of the deterministic nature of $G_n$, the desired result follows by applying the second inequality of Theorem \ref{Gronw} i.e.
	\begin{equation}
	\label{app}	\left\|\sup _{0 \leq j \leq N_h}\left|Y^j\right|\right\|_{2 p}=\left[  E\left(\sup_{0\leq j\leq N_h} |Y^j|\right)^{2p}         \right]^{1/2p}\leq \left(1+ \frac{1}{1-p}\right)^{1/2p} \left( \prod_{j=0}^{N_h-1}(1+G_j)^p\right)^{1/2p}\left(E\sup_{0\leq j\leq N_h}F_j\right)^{1/2}.
	\end{equation}
	We study first the term
	$$I = \left( \prod_{j=0}^{N_h-1}(1+G_j)^p\right)^{1/2p}.$$
	Recall that $hN_h \leq T$ and that $1+x \leq e^x$. Then,
	\begin{align*}
	&\prod_{k=0}^{N_h-1}(1+G_k)^p \leq \exp\left (p\sum_{k=0}^{N_h-1}G_k  \right) \leq \exp\left( p\sum_{k=0}^{N_h-1}2hL(1-2h_0L)^{-1}                 \right)=\exp\left( 2phN_hL(1-2h_0L)^{-1}                 \right)\\ 
	&\leq \exp\left( 2pTL(1-2h_0L)^{-1}                 \right)
	\end{align*}
	which proves that 
	\begin{equation}
	\label{expterm} I \leq \exp\left(\frac{TL }{1-2h_0L}                \right).
	\end{equation}
	For the term 
	\[
	J = \left(E\sup_{0\leq j\leq N_h}F_j\right)^{1/2}
	\]
	we can easily obtain the inequality
	\begin{equation}
	\label{const}J\leq( |Y^0|^2+2(1-2h_0L)^{-1}LT)^{1/2}.
	\end{equation}
	The result follows by combining \eqref{app}, \eqref{expterm} and \eqref{const}.
\end{proof}

\begin{remark}
	Note that the bound provided is independent of the step size $h$.
\end{remark}

\medskip

\medskip

\noindent\textbf{Acknowledgment}
	Work of the second author was supported by the Indian National Science Academy (INSA) under the ``INSA Honorary Scientist" scheme at the CR Rao Advanced Institute of Mathematics, Statistics and Computer Science, Hyderabad, India.

\end{document}